\documentclass[12pt,reqno,a4paper]{amsart}
\usepackage{amsmath,amssymb,amsfonts,amsthm,amscd}
\usepackage[mathscr]{eucal}
\usepackage{hyperref}

\author{Theodore~Th.~Voronov}
\address{School of Mathematics, University of Manchester, Sackville Street, Manchester, M60 1QD, United Kingdom}
\email{theodore.voronov@manchester.ac.uk}

\title[Mackenzie theory and $Q$-manifolds]{Mackenzie theory and $Q$-manifolds}

\newtheorem{thm}{Theorem}

\newtheorem*{genprinc}{General principle}

\newtheorem{prop}{Proposition}[section]

\theoremstyle{definition}
\newtheorem{de}{Definition}
\newtheorem{ex}{Example}[section]
\newtheorem{rem}{Remark}[section]

\def\co{\colon\thinspace}

\renewcommand{\leq}{\leqslant}

 \DeclareMathOperator{\Vect}{Vect}

\newcommand{\der}[2]{{\frac{\partial {#1}}{\partial {#2}}}}

\newcommand{\Z}{{\mathbb Z_{2}}}
\newcommand{\ZZ}{{\mathbb Z}}
\newcommand{\p}{\partial}
\newcommand{\fun}{C^{\infty}}

\newcommand{\w}{{\boldsymbol{w}}}

\renewcommand{\a}{\alpha}
\renewcommand{\b}{\beta}
\def\e{\varepsilon}

\newcommand{\g}{{\gamma}}
\newcommand{\h}{\eta}
\renewcommand{\t}{\theta}

\newcommand{\lam}{{\lambda}}
\newcommand{\x}{{\xi}}
\def\d{\delta}
\def\t{\theta}

\newcommand{\jtt}{{\tilde \jmath}}
\newcommand{\itt}{{\tilde \imath}}

\newcommand{\ut}{{\tilde u}}
\newcommand{\vt}{{\tilde v}}

\begin{document}
\begin{abstract} We give a simple characterization of Mackenzie's
double Lie algebroids in terms of homological vector fields.
Application to the `Drinfeld double' of Lie bialgebroids is given
and an extension to the multiple case is suggested.
\end{abstract}
\maketitle 

\section*{Introduction}

Double Lie algebroids arose in the works on double Lie
groupoids~\cite{mackenzie:secondorder1, mackenzie:secondorder2} and
in connection with an analog for Lie bialgebroids of the Drinfeld
double of Lie bialgebras~\cite{mackenzie:doublealg,
mackenzie:doublealg2, mackenzie:drinfeld}. They originally appeared
as the tangent objects for double Lie groupoids, and later their
properties were axiomatized to give the abstract notion. Their
immediate application, as well as that of double Lie groupoids, is
in Poisson geometry or more generally in the larger subject that may
be called \textbf{Bracket Geometry}, which embraces topics from
(multiple) groupoid and algebroid theory through homotopy algebras
to geometrical structures arising in deformation and quantization
theory.

The richness of the theory of double Lie algebroids, due to Kirill
Mackenzie,  can be seen in  numerous non-obvious structures,
isomorphisms and dualities arising in it. Notice, for example, a
non-trivial duality theory for double and triple vector bundles,
where interesting discrete symmetry groups
appear~\cite{mackenzie:duality}.

For a long time, the application of double Lie algebroids was
somewhat hindered by the complexity of their original definition. An
obvious part of the definition is, of course, the structure of two
Lie algebroids (in fact, four, on the four sides of a double vector
bundle; but they can be reduced to the two `main' ones). The
difficulty was to state a compatibility condition for them. A system
of conditions, highly non-trivial in formulation, was found by
Mackenzie~\cite{mackenzie:doublealg} as an abstraction of the Lie
functor of double Lie groupoids, and was proved to be the correct
one, for example, by showing that it is satisfied by the so-called
`cotangent doubles' of Lie bialgebroids.

In this paper we analyze Mackenzie's conditions and prove that they
are equivalent to a simple commutativity condition for homological
vector fields on a supermanifold naturally associated with a given
double vector bundle. This radically simplifies the theory and opens
ways to an immediate extension to the multiple case, i.e., $n$-fold
Lie algebroids, as well as $n$-fold `bi-' Lie algebroids. (The
latter will be the subject of a forthcoming paper with Kirill
Mackenzie~\cite{mackenzie:bidouble}.)

Our main statement (Theorem~\ref{thm.main} below) establishes
equivalence between two notions: double Lie algebroids in the sense
of Mackenzie and \textit{double Lie antialgebroids}  as defined in
this paper.

In the course of a proof, we show that Mackenzie's `Condition III'
(see Section~\ref{sec.doublealg} below), which pertains to a certain
bialgebroid, actually subsumes his other conditions.

It all fits into a \textbf{big picture}, which is as follows. For a
given double vector bundle we consider all its \textit{neighbors},
that is, the double vector bundles obtained by dualization and
parity reversion. (There are twelve of them, including the original
bundle.) We can say that a structure such as that of a double Lie
algebroid is manifested in various ways in particular structures on
all of these neighbor double vector bundles. This extends the idea
that, say, a Lie algebra  $\mathfrak g$ has equivalent
manifestations as a linear Poisson bracket (on the coalgebra
$\mathfrak g^*$), as a linear Schouten bracket (on the anticoalgebra
$\Pi\mathfrak  g^*$) and as a quadratic homological vector field (on
the antialgebra $\Pi\mathfrak g$). See, for
example,~\cite{tv:graded}. For double Lie algebroids, out of the
twelve neighbors, five allow structures with easily formulated
compatibility conditions. It has turned out that four of them are
reformulations of Mackenzie's Condition III, and the remaining one
is precisely our commutativity condition.

We wish to emphasize that in our work, supermanifolds provide
powerful tools that we apply to ordinary (``purely even'') objects.
Although we show that everything works also in a `superized'
context, this was not the main goal.

Some parts of the proofs are calculations in coordinates. They can
no doubt be replaced by coordinate-free arguments, by extending
methods used by Mackenzie.  However, I wish to note that a
`motivated' calculation in coordinates is sometimes the quickest way
to get to the crux of the matter and allows to notice facts
sometimes obscured by a more `abstract' presentation.

It has been of considerable interest among  experts to give an
alternative simpler description of double Lie algebroids since the
notion first appeared. I have always believed that such a
description should be  in terms of supermanifolds and homological
vector fields. My own earliest notes on the problem, motivated by
numerous inspiring discussions with Mackenzie, date back to 2002.
Unfortunately, this work was interrupted, which prevented me from
giving a solution at that time. I came back to it in 2003 and again
in June 2006 right before the Bia{\l}owie\.{z}a conference (see
below), wishing to discuss the problem  with Kirill Mackenzie there,
and at this time solved it completely. I wrote about the solution to
Kirill Mackenzie, Yvette Kosmann-Schwarzbach, Alan Weinstein and
Dmitry Roytenberg. Roytenberg, after learning about the statement of
Theorem~\ref{thm.main} below, told me that it was known to him, but
he did not possess a proof. As I learned from Alan Weinstein's email
(even before my work was completed), his former Ph.D. students
A.~Gracia-Saz and R.~A.~Mehta were also working on the problem;
there is a reference to a work in progress in Mehta's
thesis~\cite{mehta:thesis}, but I am not aware of any outcome.

The paper is organized as follows.

In Section~\ref{sec.doublealg} we recall the definition of double
Lie algebroids.

In Section~\ref{sec.background} we recall the description of
(ordinary) Lie algebroids in the language of homological vector
fields, and revise double vector bundles. In particular, we
introduce partial reversions of parity.

In Section~\ref{sec.main} we define double Lie antialgebroids and
give our main statement (Theorem~\ref{thm.main}).

In Section~\ref{sec.analysis} we analyze the three conditions
appearing in the definition of double Lie algebroids and give  a
proof of Theorem~\ref{thm.main}.

In Section~\ref{sec.big} we show how the equivalence of Mackenzie's
notion of double Lie algebroids and our notion of double Lie
antialgebroids is a part of a bigger picture. \textit{Modulo} some
facts established in Section~\ref{sec.analysis}, this provides an
alternative proof of Theorem~\ref{thm.main}.

In Section~\ref{sec.appl} we show the equivalence of Mackenzie's and
Roytenberg's doubles of Lie bialgebroids and discuss an extension of
the whole theory to the multiple   case.

\subsection*{Terminology and notation.}

We use the standard language of supermanifolds. The letter $\Pi$
denotes the parity reversion functor, and notation such as
$\Phi^{\Pi}$ is used for linear maps induced on the
\textit{opposite} (parity reversed) objects \footnote{For
specialists we may note that we do not make a distinction in
notation between, say, $\Pi E$ and $E\Pi$, though, practically, we
use $E\Pi$ and $\Pi E^*$ for dual vector bundles $E$ and $E^*$ to
avoid extra signs. We use left coordinates on $E$ and right
coordinates on $E^*$.}. Commutators and similar notions are always
understood in the $\Z$-graded sense. A tilde over an object is used
to denote its parity. A \textit{$Q$-manifold} means a supermanifold
endowed with a homological vector field; likewise, $P$- and
\textit{$S$-manifolds} mean those with a Poisson or Schouten (= odd
Poisson) bracket. A \textit{$QS$-manifold} means one with $Q$- and
$S$- structures that are compatible (the vector field is a
derivation of the bracket, cf.~\cite{yvette:exact}). In general,
notation and terminology are close to our paper~\cite{tv:graded}. We
wish to draw the reader's special attention to  our normally
dropping the prefix `super-' when this cannot cause   confusion and
speaking, as a rule, of `manifolds' meaning supermanifolds, `Lie
algebras' meaning superalgebras, etc.

\subsection*{Acknowledgements.}

I wish to thank my good friends Kirill Mackenzie, Yvette
Kosmann-Schwarzbach and Hovhannes Khudaverdian, for inspiring
discussions, most valuable criticism, and advice. Kirill Mackenzie
has pioneered the whole subject of multiple and bi- structures in
the groupoid and algebroid world. I  thank most cordially Yvette
Kosmann-Schwarzbach for numerous comments and remarks that helped to
improve the original manuscript. My special thanks go to the
organizers of the annual international Workshops on Geometric
Methods in Physics in Bia{\l}owie\.{z}a, notably to Anatol
Odzijewicz, for the highly inspiring atmosphere. Some of my notes on
the subject of this paper were made in Warsaw after the XXII
Bia{\l}owie\.{z}a Workshop, and the final result was reported at the
jubilee XXV Bia{\l}owie\.{z}a meeting in July 2006.

\section{Double Lie algebroids according to
Mackenzie}\label{sec.doublealg}

Double Lie algebroids were introduced by Mackenzie
in~\cite{mackenzie:doublealg, mackenzie:doublealg2}, see
also~\cite{mackenzie:drinfeld}, as the infinitesimal counterparts of
double Lie groupoids. The latter notion is a double object in the
sense of Ehresmann, i.e., a groupoid object in the category of
groupoids. Therefore, it has a natural categorical formulation.
Compared to it, the abstract notion of a double Lie algebroid is
rather complicated and non-obvious. One of the reasons for this, is
that properties of brackets for Lie algebroids are not expressed
diagrammatically, so one cannot approach double objects for them by
methods of category theory. Mackenzie's conditions (see below) come
about as an abstraction of the properties of the double Lie
algebroid of a double Lie groupoid discovered
in~\cite{mackenzie:secondorder1, mackenzie:secondorder2}.

\begin{de}
\label{def.dlie} A double vector bundle
\begin{equation} \label{eq.dvb}
    \begin{CD} D@>>> B\\
                @VVV  @VVV \\
                A@>>>M
    \end{CD}
\end{equation}
is a \textit{double Lie algebroid} if all sides are (ordinary) Lie
algebroids and the following \textbf{conditions I}, \textbf{II}, and
\textbf{III} are satisfied:
\begin{description}
  \item[Condition I] With  respect to the vertical structures of
  Lie algebroids, $\begin{array}{c} D\\
                \downarrow  \\
                A
    \end{array}$ and
    $\begin{array}{c} B\\
                \downarrow  \\
                M
    \end{array}$, all maps related with the horizontal vector bundle
    structures are Lie algebroid morphisms (more precisely, it includes the projections, the zero sections,
    the fiberwise addition and multiplication by scalars). The same
    holds with vertical/horizontal structures interchanged.
  \item[Condition II] The horizontal arrows in the diagram
  \begin{equation*}
    \begin{CD} D@>a>> TB\\
                @VVV  @VVV \\
                A@>a>>TM
    \end{CD}
\end{equation*}
where at the right there is the tangent prolongation of the Lie
algebroid $B\to M$, and $a$ stands for the anchors, define a Lie
algebroid morphism. The same holds with vertical/horizontal
structures  interchanged.
  \item[Condition III] The vertical arrows in the diagram
  \begin{equation*}
    \begin{CD} D^{*A}@>>> K^*\\
                @VVV  @VVV \\
                A@>>>M
    \end{CD}
\end{equation*}
define a Lie algebroid morphism. Here $K$ denotes the core. The same
holds with vertical/horizontal, and $A$ and $B$,  interchanged. The
vector bundles in duality $D^{*A}\to K^*$ and $D^{*B}\to K^*$ define
a Lie bialgebroid.
\end{description}
(An explication of these conditions will be given below.)
\end{de}

Recall that a \textit{double vector bundle} such as~\eqref{eq.dvb}
is defined by the condition that all vector bundle structure maps in
one direction (horizontal or vertical) are vector bundle morphisms
for another direction. The \textit{core} $K$  is defined as the
intersection of the kernels of the projections $D\to A$ and $D\to B$
considered as vector bundle morphisms (w.r.t.  the other structure).
$K$ is a vector bundle over $M$. It is a theorem due to Mackenzie
that taking the two duals of $D$ considered as a vector bundle
either over $A$ or over $B$ leads to two double vector bundles
\begin{equation*}
    \begin{CD} D^{*A}@>>> K^*\\
                @VVV  @VVV \\
                A@>>>M
    \end{CD} \text{\quad and \quad }
    \begin{CD} D^{*B}@>>> B\\
                @VVV  @VVV \\
                K^*@>>>M
    \end{CD}
\end{equation*}
where the vector bundles $D^{*A}\to K^*$ and $D^{*B}\to K^*$ over
the co-core $K^*$ are\,---\,unexpectedly\,---\,in a natural duality.
All these facts, as well as the notion of the \textit{tangent
prolongation} of a Lie algebroid, can be found
in~\cite[Ch.~9]{mackenzie:book2005}, see also
\cite{mackenzie:duality, mackenzie:bialg}. \textit{Lie bialgebroids}
were introduced by Mackenzie and Xu~\cite{mackenzie:bialg}. Their
theory was  advanced by Y.~Kosmann-Schwarzbach~\cite{yvette:exact},
who in particular gave a very handy form of the definition, which we
use. See~\cite{mackenzie:book2005}.

\section{Lie algebroids and double vector bundles: some background}
\label{sec.background}

In this section we  develop tools that will be later used for an
alternative description of double Lie algebroids (our main goal).

Henceforth we work in the `super' setup, i.e., we consider
supermanifolds and bundles of supermanifolds. However, we
systematically skip the prefix `super-' except when we wish to make
an emphasis.   All the constructions from the previous section carry
over to the super case.

We  use \textit{graded manifolds} as defined in~\cite{tv:graded},
i.e., supermanifolds endowed with an extra $\ZZ$-grading in the
algebras of functions, in general not related with parity. We refer
to such grading as \textit{weight}.

Let us recall some known facts concerning \textbf{Lie algebroids}.

It was first  shown by Vaintrob~\cite{vaintrob:algebroids} that Lie
algebroids can be described by homological vector fields. We shall
recall this correspondence using the description given
in~\cite{tv:graded} in the language of derived brackets. As
mentioned, we consider the `superized' version (i.e., `super' Lie
algebroids) by default.

Let $F\to M$ be a vector bundle. The total space $F$ is naturally a
graded manifold,  the (pullbacks of) functions on the base $M$
having weight $0$ and linear functions on the fibers, weight $1$.
Using weights is very helpful for describing various geometric
objects. For example, vector fields  of weight $-1$ on $F$
correspond to sections of $F$ (or $\Pi F$, see below). Vector fields
of weight $0$ are generators of fiberwise linear transformations.
Vector fields of weight $1$ can be used to generate brackets of
sections. More precisely:  a \textit{Lie antialgebroid} structure on
$F\to M$, by definition,  is given by a homological field $Q\in
\Vect (F)$ of weight $1$.

\textit{There is a one-to-one correspondence between Lie
antialgebroids and Lie algebroids}, as follows.

Let $\Pi$ denote the parity reversion functor, and $F=\Pi E$ for a
vector bundle $E\to M$. Then $F$ is a Lie antialgebroid if and only
if $E$ is a Lie algebroid. The anchor and the bracket  for the
sections of $E$ are given by the following formulas:
\begin{equation}\label{eq.auf}
    a(u)f:=\bigl[[Q,i(u)],f\bigr]
\end{equation}
and
\begin{equation} \label{eq.uv}
    i([u,v]):=(-1)^{\ut}\bigl[[Q,i(u)],i(v)\bigr].
\end{equation}
Here $f\in \fun(M)$, and $u,v\in \fun(M, E)$ are sections. We use
the natural odd injection $i\co \fun(M, E)\to \Vect (\Pi E)$, which
sends a section $u\in \fun(M, E)$ to a vector field $i(u)\in \Vect
(\Pi E)$   of weight $-1$. The map $i$ is an odd isomorphism between
the space of sections $\fun(M, E)$ and the subspace $\Vect_{-1}(\Pi
E)\subset \Vect(\Pi E)$ of all vector fields of weight $-1$. By
counting weights, one can see that the LHS's of~\eqref{eq.auf}
and~\eqref{eq.uv} are well-defined. The properties of the bracket
and anchor are deduced from the identity $Q^2=0$ as    is standard
in the derived brackets method. Conversely, starting from a Lie
algebroid structure in $E\to M$, one can reconstruct $Q$ on $\Pi E$
with the desired properties.

All these facts can be checked without coordinates; however,
introducing local coordinates makes them particularly transparent.
Let $x^a$ denote local coordinates on the base $M$. We shall use
$u^i$ and $\x^i$ for linear coordinates in the fibers of $E$ and
$F=\Pi E$, respectively. Changes of coordinates have the following
form:
\begin{align*}
    x^a&=x^a(x'),\\
    u^i&=u^{i'}T_{i'}{}^{i}(x'),
    \intertext{and}
\x^i&=\x^{i'}T_{i'}{}^{i}(x').
\end{align*}
The map $i\co \fun(M, E)\to \Vect (\Pi E)$ has the following
appearance in coordinates\footnote{Here the components $u^i(x)$ of a
section $u=u^i(x)e_i$ should not be confused with the coordinates
$u^i$, which are functions on the total space $E$. In particular,
each of the coordinates $u^i$ has a certain fixed parity, and the
correspondent  component $u^i(x)$ of a section $u$ has the same or
the opposite parity as the coordinate $u^i$ depending on whether the
section $u$ is even or odd.}:  $i\co u=u^i(x)e_i\mapsto i(u)$, where
\begin{equation}\label{eq.iu}
    i(u)=(-1)^{\ut}u^i(x)\der{}{\x^i}\,.
\end{equation}
Clearly, the RHS of \eqref{eq.iu} is the general form of a vector
field of weight $-1$ on $F$. A vector field $Q$ of weight $1$ on $F$
in coordinates has the form
\begin{equation*}
    Q=\x^iQ_i^a(x)\,\der{}{x^a}+\frac{1}{2}\,\x^i\x^j Q_{ji}^k
    (x)\,\der{}{\x^k}\,.
\end{equation*}
Equations~\eqref{eq.auf} and \eqref{eq.uv} produce the following
formulas for the anchor:
\begin{equation*}
    a(u)=u^i(x)\,Q_i^a(x)\,\der{}{x^a}\,,
\end{equation*}
and for the brackets:
\begin{equation*}
    [u,v]=\Bigl(u^i Q_i^a\,\p_av^k-(-1)^{\ut(\vt+1)}v^i Q_i^a\,\p_au^k-
    (-1)^{\itt(\vt+1)} u^iv^j Q_{ji}^k\Bigr)e_k\,,
\end{equation*}
where we abbreviated $\p_a=\p/\p x^a$. In particular, for the
elements of the local frame $e_i$ we have
\begin{equation*}
    [e_i,e_j]=(-1)^{\jtt} Q_{ij}^k(x)\,e_k\,.
\end{equation*}

Let us now proceed to \textbf{double vector bundles}. Let
\begin{equation} \label{eq.dvb2}
    \begin{CD} D@>>> B\\
                @VVV  @VVV \\
                A@>>>M
    \end{CD}
\end{equation}
be a double vector bundle (in the category of supermanifolds). The
manifold $D$ is naturally bi-graded, by weights corresponding to the
two vector bundle structures. If necessary we denote these weights
by $\w_1$ and $\w_2$, or by $\w_A$ and $\w_B$, as convenient.

Double vector bundles allow fiberwise reversion of parity in both
directions, horizontal and vertical. We denote the corresponding
operations by $\Pi_1$ and $\Pi_2$ (or by $\Pi_A$ and $\Pi_B$ when
convenient). (Such operations should be studied together with the
operations of dualization in Mackenzie's
theory~\cite{mackenzie:duality}.) For example, for a double vector
bundle  given by~\eqref{eq.dvb2}, the vertical reversion of parity
$\Pi_1=\Pi_A$ gives
\begin{equation}\label{eq.pa}
    \begin{CD} \Pi_A D@>>> \Pi B\\
                @VVV  @VVV \\
                A@>>>M
    \end{CD}
\end{equation}
which is a new double vector bundle. One can apply horizontal
reversion of parity to it (applying the vertical reversion again
takes us back), or do it  the other way round.

\begin{prop} \label{prop.pp} The operations $\Pi_1$ and $\Pi_2$ commute:
\begin{equation*}
    \Pi_1\Pi_2=\Pi_2\Pi_1.
\end{equation*}
More precisely, for a double vector bundle given by~\eqref{eq.dvb2},
there is a natural isomorphism of double vector bundles
\begin{equation} \label{eq.pi1pi2}
    \begin{CD} \Pi_B\Pi_A D@>>> \Pi B\\
                @VVV  @VVV \\
                \Pi A@>>>M
    \end{CD}
    \text{\quad and \quad}
    \begin{CD} \Pi_A\Pi_B D@>>> \Pi B\\
                @VVV  @VVV \\
                \Pi A@>>>M
    \end{CD}
\end{equation}
where we used the more suggestive notation $\Pi_A=\Pi_1$ and
$\Pi_B=\Pi_2$.
\end{prop}

Denote the common value of the ultimate total spaces
in~\eqref{eq.pi1pi2} by $\Pi^2 D$,
\begin{equation*}
    \Pi^2 D:= \Pi_A\Pi_B D = \Pi_B\Pi_A D\,,
\end{equation*}
and call it the \textit{complete reversion of parity} of $D$.

\begin{rem} Everything extends immediately to the  $n$-fold
case, with $\Pi^n D$ being the \textit{complete parity reversion} of
the ultimate total space $D$ of an $n$-fold vector bundle. There are
partial parity reversion operations $\Pi_{r}$ such that
$\Pi_{r}\Pi_{s}=\Pi_{s}\Pi_{r}$ and $\Pi^n D=\Pi_1\cdot \ldots
\cdot\Pi_n D$.  (See Section~\ref{sec.appl}.)
\end{rem}

A coordinate language, particularly handy for visualizing double
(and multiple) vector bundles, is as follows. Consider a double
vector bundle given by~\eqref{eq.dvb2}. As above, denote local
coordinates on $M$ by $x^a$. Let $u^i$ and $w^{\a}$ be linear
coordinates on the fibers of $A\to M$ and $B\to M$, respectively. On
$D$ we have coordinates $x^a, u^i, w^{\a}, z^{\mu}$ so that $u^i,
z^{\mu}$ are linear fiber coordinates for $D\to B$ and $w^{\a},
z^{\mu}$, for $D\to A$. Coordinate changes have the form:
\begin{align}
    x^a&=x^a(x'),\\
    u^i&=u^{i'}T_{i'}{}^{i}(x'),\label{eq.lawforu}\\
    w^{\a}&=w^{\a'}T_{\a'}{}^{\a}(x'),\\
    z^{\mu}&=z^{\mu'}T_{\mu'}{}^{\mu}(x')+u^{i'}w^{\a'}T_{\a'
    i'}{}^{\mu}(x')\,.\label{eq.lawforz}
\end{align}
This is a convenient description of a double vector bundle
structure. (The reader who lacks a taste for coordinate calculations
may  translate it into a language of local trivializations.) In
particular, for two weights we have $\w_1=\# u+\# z$ and $\w_2=\#
w+\# z$, where $\#$ denotes the degree in the respective variable.
Everything extends directly to multiple vector bundles.

Partial parity reversion is described as follows. For~\eqref{eq.pa}
we have $x^a, u^i, \h^{\a}, \t^{\mu}$ as coordinates on $\Pi_AD$, so
that $\widetilde{\h^{\a}}=\widetilde{w^{\a}}+1=\widetilde{\a}+1$,
$\widetilde{\t^{\mu}}=\widetilde{z^{\mu}}+1=\widetilde{\mu}+1$, and
changes of coordinates have the form
\begin{align*}
    u^i&=u^{i'}T_{i'}{}^{i},\\
    \h^{\a}&=\h^{\a'}T_{\a'}{}^{\a},\\
    \t^{\mu}&=\t^{\mu'}T_{\mu'}{}^{\mu}+(-1)^{\itt'}u^{i'}\h^{\a'}T_{\a'
    i'}{}^{\mu}\,,
\end{align*}
where we suppress coordinates on $M$. Here $\h^{\a}, \t^{\mu}$ are
fiber coordinates for $\Pi_A D\to A$ and $\h^{\a}$ are fiber
coordinates for $\Pi B\to M$. Similarly, for
\begin{equation}\label{eq.pb}
    \begin{CD} \Pi_B D@>>>  B\\
                @VVV  @VVV \\
                \Pi A@>>>M
    \end{CD}
\end{equation}
we have $x^a, \x^i, w^{\a}, \t^{\mu}$ as coordinates on $\Pi_BD$,
where $\widetilde{\x^{i}}=\widetilde{u^{i}}+1={\itt}+1$ and
$\widetilde{\t^{\mu}}=\widetilde{z^{\mu}}+1=\widetilde{\mu}+1$, with
changes of coordinates of the form
\begin{align*}
    \x^i&=\x^{i'}T_{i'}{}^{i},\\
    w^{\a}&=w^{\a'}T_{\a'}{}^{\a},\\
    \t^{\mu}&=\t^{\mu'}T_{\mu'}{}^{\mu}+\x^{i'}w^{\a'}T_{\a'
    i'}{}^{\mu}\,.
\end{align*}
(Note different meanings of $\t^{\mu}$ for $\Pi_A D$ and $\Pi_B D$.)
Applying parity reversion once again we obtain
\begin{equation} \label{eq.pibpia}
    \begin{CD} \Pi_B\Pi_A D@>>> \Pi B\\
                @VVV  @VVV \\
                \Pi A@>>>M
    \end{CD}
\end{equation}
with coordinates on $x^a, \x^i, \h^{\a}, z^{\mu}$ on $\Pi_B\Pi_A D$
with the transformation law
\begin{align*}
    \x^i&=\x^{i'}T_{i'}{}^{i},\\
    \h^{\a}&=\h^{\a'}T_{\a'}{}^{\a},\\
    z^{\mu}&=z^{\mu'}T_{\mu'}{}^{\mu}+(-1)^{\itt'}\x^{i'}\h^{\a'}T_{\a'
    i'}{}^{\mu}\,,
\end{align*} and
\begin{equation} \label{eq.piapib}
    \begin{CD} \Pi_A\Pi_B D@>>> \Pi B\\
                @VVV  @VVV \\
                \Pi A@>>>M
    \end{CD}
\end{equation}
with coordinates $x^a, \x^i, \h^{\a}, z^{\mu}$ on $\Pi_A\Pi_B D$
with the transformation law
\begin{align*}
    \x^i&=\x^{i'}T_{i'}{}^{i},\\
    \h^{\a}&=\h^{\a'}T_{\a'}{}^{\a},\\
    z^{\mu}&=z^{\mu'}T_{\mu'}{}^{\mu}+(-1)^{\itt'+1}\x^{i'}\h^{\a'}T_{\a'
    i'}{}^{\mu}\,.
\end{align*}
Note that these transformation laws are the same up to a change of
sign of the $z$-coordinate. In particular, this gives a proof for
Proposition~\ref{prop.pp}.

\section{Main statement}\label{sec.main}

In this section we give our main statement, which is a
characterization of Mackenzie's double Lie algebroids in terms of
graded $Q$-manifolds.

\begin{de} \label{def.dalie} A double vector bundle
\begin{equation} \label{eq.dvb3}
    \begin{CD} H@>>> G\\
                @VVV  @VVV \\
                F@>>>M
    \end{CD}
\end{equation}
is a \textit{double Lie antialgebroid} if it is endowed with two
homological vector fields $Q_1$ and $Q_2$ on the manifold $H$ of
weights $(1,0)$ and $(0,1)$, respectively, such that
\begin{equation}\label{eq.q1q2}
    [Q_1,Q_2]=0.
\end{equation}
\end{de}

\begin{rem}
Equivalently, a double Lie antialgebroid can be defined as a double
vector bundle such as~\eqref{eq.dvb3} with a (single) homological
vector field $Q$ on $H$ of \textbf{total} weight $1$. Then $Q$ can
be decomposed into the sum $Q_1+Q_2$ of  fields of weights $(1,0)$
and $(0,1)$, and $Q^2=0$ is equivalent to $Q_1^2=Q_2^2=[Q_1,
Q_2]=0$.
\end{rem}

\begin{rem} \label{rem.mantialg}
An extension to \textit{multiple Lie antialgebroids} is immediate.
An \textit{$n$-fold Lie antialgebroid} is an $n$-fold vector bundle,
which therefore gives rise to an $n$-graded structure on its
(ultimate) total space, endowed with a homological vector field $Q$
of total weight $1$. In  greater detail, we obtain $n$ commuting
homological fields $Q_{r}$, $r=1,\ldots,n$, on the ultimate total
space of weights $(0,\ldots, 1, \ldots, 0)$, respectively. (Here $1$
stands at the $r$-th place, all other weights being zero.) Then
\begin{equation*}
    Q=Q_1+\ldots+Q_n.
\end{equation*}
We shall come back to this in Section~\ref{sec.appl}.
\end{rem}

\begin{ex} Let $(E, E^*)$ be a Lie
bialgebroid with base $M$. Then
\begin{equation} \label{eq.dimadouble}
    \begin{CD} T^*\Pi E=T^*\Pi E^*@>>> \Pi E^*\\
                @VVV  @VVV \\
                \Pi E@>>>M
    \end{CD}
\end{equation}
is a double vector bundle, as one can check. This is a superization
of Mackenzie~\cite{mackenzie:doublealg, mackenzie:doublealg2}. The
natural diffeomorphism in the upper-left corner
of~\eqref{eq.dimadouble} is a superized version of a theorem of
Mackenzie and Xu~\cite{mackenzie:bialg} extending a statement of
Tulczyjew~\cite{tulczyjew:1977}. The Lie algebroid structures in $E$
and $E^*$ give rise to homological fields on $\Pi E$ and $\Pi E^*$,
respectively.  These two vector fields correspond to two functions
(`linear Hamiltonians') on $T^*\Pi E=T^*\Pi E^*$, of weights $(1,0)$
and $(0,1)$. That $(E,E^*)$ is a  bialgebroid  is equivalent to the
commutativity of these Hamiltonians (due to
Roytenberg~\cite{roytenberg:thesis}, see also \cite{tv:graded}).
Therefore the corresponding Hamiltonian vector fields make the
double vector bundle~\eqref{eq.dimadouble}  a double Lie
antialgebroid. We shall come back to this example in
Section~\ref{sec.appl}.
\end{ex}

\begin{thm}[A Characterization of Double Lie
Algebroids]\label{thm.main}
A double Lie algebroid structure  in a double vector bundle such
as~\eqref{eq.dvb2} is equivalent to a double Lie antialgebroid
structure in the complete parity reversion double vector bundle
\begin{equation} \label{eq.pikvadrat}
    \begin{CD} \Pi^2 D@>>> \Pi B\\
                @VVV  @VVV \\
                \Pi A@>>>M
    \end{CD}
\end{equation}
i.e., to a homological field $Q=Q_1+Q_2$ on $\Pi^2 D$ of total
weight $1$.
\end{thm}

To appreciate the statement one may compare the three conditions of
Definition~\ref{def.dlie} with the commutativity
equation~\eqref{eq.q1q2} of Definition~\ref{def.dalie}.

Let us show how a homological vector field $Q=Q_1+Q_2$ of total
weight $1$ on $\Pi^2 D$ generates  Lie algebroid structures on all
sides of~\eqref{eq.dvb2}. As in our discussion of a (single) Lie
algebroid above, everything can be formulated in a coordinate-free
setting. However, using coordinates sheds some extra light.

For the sake of concreteness, we consider $\Pi^2 D=\Pi_B\Pi_A D$
with natural coordinates $x^a, \x^i, \h^{\a}, z^{\mu}$ thereon. A
vector field $Q_1\in \Vect (\Pi^2 D)$ of weight $(1,0)$ has the form
\begin{multline}\label{eq.fieldqa}
    Q_1=\x^iQ_i^a\der{}{x^a}+\frac{1}{2}\x^i\x^jQ_{ji}^k\der{}{\x^k} \\
    +\left(\x^i\h^{\a}Q_{\a i}^{\b}+z^{\mu}Q_{\mu}^{\b}\right)\der{}{\h^{\b}}
    +\left(\frac{1}{2}\x^i\x^j\h^{\a}Q_{\a ji}^{\lambda}+\x^iz^{\mu}Q_{\mu i}^{\lambda}\right)
    \der{}{z^{\lambda}}\,,
\end{multline}
while $Q_2\in \Vect (\Pi^2 D)$ of weight $(0,1)$, the form
\begin{multline}\label{eq.fieldqb}
    Q_2=\h^{\a}Q_{\a}^a\der{}{x^a}+
    \left(\h^{\a}\x^iQ_{i\a}^{j}+z^{\mu}Q_{\mu}^{j}\right)\der{}{\x^{j}}\\
    + \frac{1}{2}\h^{\a}\h^{\b}Q_{\b\a}^{\g}\der{}{\h^{\g}}+
    \left(\frac{1}{2}\h^{\a}\h^{\b}\x^{i}Q_{i \b\a}^{\lambda}+\h^{\a}z^{\mu}Q_{\mu\a}^{\lambda}\right)
    \der{}{z^{\lambda}}\,.
\end{multline}
All coefficients here are functions of $x^a$. Now, \textit{due to
the fact that $Q_1$ has weight $0$ w.r.t. the horizontal fiber
coordinates, it admits partial parity reversion in this direction},
giving a vector field on $\Pi_B D$:
\begin{multline}\label{eq.fieldqapi}
    Q_1^{\Pi}:=\x^iQ_i^a\der{}{x^a}+\frac{1}{2}\x^i\x^jQ_{ji}^k\der{}{\x^k} \\
    +\left(\x^i w^{\a}Q_{\a i}^{\b}+\t^{\mu}Q_{\mu}^{\b}\right)\der{}{w^{\b}}
    +\left(\frac{1}{2}\x^i\x^j w^{\a}Q_{\a ji}^{\lambda}+\x^i\t^{\mu}Q_{\mu\a}^{\lambda}\right)
    \der{}{\t^{\lambda}}\,.
\end{multline}
Similarly, $Q_2$ \textit{allows vertical parity reversion}, which
gives
\begin{multline}\label{eq.fieldqbpi}
    Q_2^{\Pi}:=\h^{\a}Q_{\a}^a\der{}{x^a}+
    \left(\h^{\a}u^iQ_{i\a}^{j}+\t^{\mu}Q_{\mu}^{j}\right)\der{}{u^{j}}\\
    + \frac{1}{2}\h^{\a}\h^{\b}Q_{\b\a}^{\g}\der{}{\h^{\g}}+
    \left(\frac{1}{2}\h^{\a}\h^{\b}u^{i}Q_{i \b\a}^{\lambda}+\h^{\a}\t^{\mu}Q_{\mu\a}^{\lambda}\right)
    \der{}{\t^{\lambda}}
\end{multline}
on $\Pi_A D$.

Both $Q_1^{\Pi}$ and $Q_2^{\Pi}$ are homological fields. They define
Lie antialgebroid structures on the vector bundles $\Pi_BD\to B$ and
$\Pi_A D\to A$, which correspond to Lie algebroid structures on
$D\to B$ and $D\to A$.

The restrictions of $Q_1^{\Pi}$ and $Q_2^{\Pi}$ on $\Pi A$ and $\Pi
B$, respectively, treated as submanifolds (zero sections) in
$\Pi_BD$ and $\Pi_A D$ are tangent to these submanifolds and define
homological vector fields
\begin{equation}\label{eq.fieldqapi2}
    Q_1^{(0)}=\x^iQ_i^a\der{}{x^a}+\frac{1}{2}\x^i\x^jQ_{ji}^k\der{}{\x^k}
\end{equation}
on $\Pi A$ and
\begin{equation}\label{eq.fieldqbpi2}
    Q_2^{(0)}=\h^{\a}Q_{\a}^a\der{}{x^a}
    + \frac{1}{2}\h^{\a}\h^{\b}Q_{\b\a}^{\g}\der{}{\h^{\g}}
\end{equation}
on $\Pi B$. This leads to Lie algebroid structures on $A\to M$ and
$B\to M$.

\section{Analysis of Mackenzie's conditions}\label{sec.analysis}

To prove Theorem~\ref{thm.main}, we can do the following. Consider a
double vector bundle given by~\eqref{eq.dvb}. Assume that all four
sides are Lie algebroids and describe them by homological vector
fields. Then we can study conditions I, II and III of
Definition~\ref{def.dlie} and see what they mean in terms of these
fields.

Recall the notion of a \textit{Lie algebroid morphism}. It is
non-obvious for the case of different bases.
See~\cite[$\S$4.3]{mackenzie:book2005} for the definition. Instead
of it, we shall use the  following  statement:

\begin{prop}\label{prop.lamorphism} Suppose $E_1\to M_1$ and $E_2\to M_2$ are
Lie algebroids defined by homological vector fields $Q_1\in \Vect \Pi E_1$ and
$Q_2\in \Vect \Pi E_2$.
A vector bundle map given by the horizontal
arrows of
\begin{equation*}
    \begin{CD} E_1 @>{\Phi}>> E_2 \\
                @VVV    @VVV\\
                M_1 @>{\varphi}>> M_2
    \end{CD}
\end{equation*}
is a morphism of Lie algebroids
if and only if  the vector  fields $Q_1$ and $Q_2$ are
$\Phi^{\Pi}$-related, where $$\Phi^{\Pi}\co \Pi E_1\to \Pi E_2$$ is
the  induced map of the opposite vector bundles.
\end{prop}

This statement first appeared, without a proof, in
Vaintrob~\cite{vaintrob:algebroids}. In our language, it is
equivalent to saying that the map $\Phi^{\Pi}\co \Pi E_1\to \Pi E_2$
is a morphism of Lie antialgebroids. It is much easier to handle
than the original definition of morphisms of Lie algebroids.

Recall that vector fields on (super)manifolds are \textit{related}
by a smooth map $F$ if the pull-back of functions $F^*$ commutes
with them. In terms of the local flows $g_t$ and $h_t$ generated by
these fields, this means that $F$ commutes with the  flows:
$g_t\circ F=F\circ h_t$.

Let us introduce a necessary notation. We have `horizontal' Lie
algebroid structures, i.e., in the vector bundles $D\to B$ and $A\to
M$, and `vertical', i.e., in $D\to A$ and $B\to M$. Hence we have
homological vector fields $Q_{DB}\in \Vect (\Pi_BD)$,  $Q_{AM}\in
\Vect (\Pi A)$,  $Q_{DA}\in \Vect (\Pi_AD)$ and $Q_{BM}\in \Vect
(\Pi B)$. If we use the notation for coordinates from
Section~\ref{sec.background}, then $Q_{DB}$ has weight $1$ in
variables $\x^i,\t^{\mu}$ and $Q_{AM}$, in variables $\x^i$.
Similarly, $Q_{DA}$ has weight $1$ in variables $\h^{\a},\t^{\mu}$
and $Q_{BM}$, in variables $\h^{\a}$.

We shall study conditions I, II and III one by one.

\subsection{Condition I}
{Condition I} is the easiest for  analysis.

Consider for concreteness the horizontal algebroid structures.
Condition I requires that all vertical structure maps: bundle
projections, zero sections, fiberwise addition and multiplication by
scalars, give morphisms of Lie algebroids. We have the following
diagrams to analyze:

\begin{equation*}
    \begin{CD} D @>>> B \\
                @V{p}VV    @VV{p}V\\
                A @>>> M
    \end{CD}
\text{\quad and \quad }
\begin{CD} D @>>> B \\
                @A{i}AA    @AA{i}A\\
                A @>>> M
    \end{CD}
\end{equation*}
for morphisms
\begin{equation*}
    \begin{CD}
    \begin{CD}D @>>>B \end{CD}\\
    @VVV \\
     \begin{CD}A @>>>M \end{CD}
     \end{CD}
    \text{\quad and \quad }
    \begin{CD}
    \begin{CD}D @>>>B \end{CD}\\
    @AAA \\
     \begin{CD}A @>>>M \end{CD}
     \end{CD}
\end{equation*}
and
\begin{equation*}
    \begin{CD} D \times_A D @>>> B\times_M B \\
                @V{+_A}VV    @VV{+_M}V\\
                D @>>> B
    \end{CD}
\text{\quad and \quad }
\begin{CD} D @>>> B \\
                @V{t_A}VV    @VV{t_M}V\\
                D @>>> B
    \end{CD}
\end{equation*}
for morphisms
\begin{equation*}
    \begin{CD}
    \begin{CD}D \times_A D @>>> B\times_M B \end{CD}\\
    @VVV \\
     \begin{CD}D @>>> B \end{CD}
     \end{CD}
    \text{\quad and \quad }
    \begin{CD}
    \begin{CD}D @>>>B \end{CD}\\
    @VVV \\
     \begin{CD}D @>>>B \end{CD}
     \end{CD}
\end{equation*}

In the language of homological vector fields, we see that the flows
generated by the vector fields $Q_{DB}$ and $Q_{AM}$ should commute
with all the vertical structure maps above, more precisely, with the
maps induced on the total spaces of the parity reversed horizontal
vector bundles. Commuting with the projection means that (the flow
of) the vector field $Q_{DB}$ acts fiberwise on the total space of
$\Pi_BD\to \Pi A$ and induces  on $\Pi A$ (the flow of) the vector
field $Q_{AM}$. Hence  $Q_{AM}$ is completely determined by
$Q_{DB}$. Consider the action of the flow of $Q_{DB}$ on the fibers
of $\Pi_BD\to \Pi A$. Commutativity with   fiberwise multiplication
by scalars, $t_{A}\co \Pi_BD \to \Pi_BD$, and addition, ${+}_{A}\co
\Pi_BD \times_{\Pi A} \Pi_BD \to \Pi_BD$, means that the flow of
$Q_{DB}$ is fiberwise linear (over $\Pi A$). This is equivalent to
the vector field $Q_{DB}$ having weight $0$ w.r.t. fiber coordinates
on $\Pi_BD\to \Pi A$. Commutativity with the zero section $\Pi A\to
\Pi_BD$ then comes about automatically.

We may summarize: if the horizontal Lie algebroid structures are
described by homological vector fields  $Q_{DB}$ and $Q_{AM}$, then
\textit{Condition I of Definition~\ref{def.dlie} is equivalent to
$Q_{DB}$ having vertical weight $0$ (its horizontal weight is $1$)
and  $Q_{AM}$ being the restriction of $Q_{DB}$ to the base $\Pi
A\subset \Pi_BD$}.

In coordinates we obtain the following expressions for $Q_{DB}$ and
$Q_{AM}$:
\begin{multline}\label{eq.qdb}
    Q_{DB}=\x^iQ_i^a\der{}{x^a}+\frac{1}{2}\x^i\x^jQ_{ji}^k\der{}{\x^k} \\
    +\left(\x^i w^{\a}Q_{\a i}^{\b}+\t^{\mu}Q_{\mu}^{\b}\right)\der{}{w^{\b}}
    +\left(\frac{1}{2}\x^i\x^j w^{\a}Q_{\a ji}^{\lambda}+\x^i\t^{\mu}Q_{\mu i}^{\lambda}\right)
    \der{}{\t^{\lambda}}\,,
\end{multline}
a vector field of weight $(1,0)$ on $\Pi_B D$, and
\begin{equation}\label{eq.qam}
    Q_{AM}=\x^iQ_i^a\der{}{x^a}+\frac{1}{2}\x^i\x^jQ_{ji}^k\der{}{\x^k}\,,
\end{equation}
a vector field on $\Pi A$ of weight $1$, the restriction of $Q_{DB}$
on $\Pi A$.

In the same way we deduce, for  the vector fields $Q_{DA} \in \Vect
(\Pi_A D)$ and $Q_{BM} \in \Vect (\Pi B)$ describing vertical Lie
algebroid structures, that \textit{ $Q_{DA}$ should have weight
$(0,1)$ on $\Pi_AD$ and $Q_{BM}$  be its restriction to $\Pi B$}.
Hence the coordinate expressions:
\begin{multline}\label{eq.qda}
    Q_{DA}=\h^{\a}Q_{\a}^a\der{}{x^a}+
    \left(\h^{\a}u^iQ_{i\a}^{j}+\t^{\mu}Q_{\mu}^{j}\right)\der{}{u^{j}}\\
    + \frac{1}{2}\h^{\a}\h^{\b}Q_{\b\a}^{\g}\der{}{\h^{\g}}+
    \left(\frac{1}{2}\h^{\a}\h^{\b}u^{i}Q_{i \b\a}^{\lambda}+\h^{\a}\t^{\mu}Q_{\mu\a}^{\lambda}\right)
    \der{}{\t^{\lambda}}\,,
\end{multline}
and
\begin{equation}\label{eq.qbm}
    Q_{BM}=\h^{\a}Q_{\a}^a\der{}{x^a}
    + \frac{1}{2}\h^{\a}\h^{\b}Q_{\b\a}^{\g}\der{}{\h^{\g}}\,.
\end{equation}

We have recovered
formulas~\eqref{eq.fieldqapi}--\eqref{eq.fieldqbpi2}.

\begin{rem}
Note that homological vector fields $Q_{DB}$ and $Q_{DA}$
determining the horizontal and vertical Lie algebroid structures are
defined on different supermanifolds, $\Pi_BD$ and $\Pi_AD$. The
crucial fact is that both have weight zero ``in the other
direction''. In particular, this allows to additionally reverse
parity and obtain vector fields defined on a common domain $\Pi^2 D$
(so that the commutativity condition  makes sense).
\end{rem}

\subsection{Condition II}

Consider Condition II of Definition~\ref{def.dlie}. For diagram
\begin{equation} \label{eq.datotbtm1}
    \begin{CD} D@>a>> TB\\
                @VVV  @VVV \\
                A@>a>>TM
    \end{CD}
\end{equation}
which is supposed to give a Lie algebroid morphism
\begin{equation}\label{eq.datotbtm}
\begin{CD}
    \begin{CD} D\\
                @VVV \\
                A
    \end{CD}
    @>>>
    \begin{CD}  TB\\
                 @VVV \\
                TM
    \end{CD}
    \end{CD}
\end{equation}
we first need to explicate the tangent prolongation  Lie algebroid
$TB\to TM$. The definition is in~\cite[$\S$9.7]{mackenzie:book2005}.
However, we shall use the following proposition instead.
\begin{prop} The tangent prolongation Lie algebroid $TE\to TM$ of a
Lie algebroid $E\to M$ is given by  the tangent prolongation of the
corresponding homological vector field $Q\in \Vect (\Pi E)$, which
is an (automatically homological) vector field on $T(\Pi E)=\Pi_{TM}
TE$.
\end{prop}

Note that for any vector bundle $E\to M$ taking tangents leads to a
double vector bundle
\begin{equation*}
    \begin{CD} TE@>>> E\\
                @VVV  @VVV \\
                TM@>>>M
    \end{CD}
\end{equation*}
(see, e.g.,~\cite[$\S$3.4]{mackenzie:book2005}), so partial parity
reversions make sense.

By differentiating the field $Q_{BM}$ given by~\eqref{eq.qbm}, we
immediately obtain a vector field $\hat Q_{BM}$ on $T(\Pi
B)=\Pi_{TM} TB$,
\begin{multline}\label{eq.qbmhat}
    \hat Q_{BM}=\h^{\a}Q_{\a}^a\der{}{x^a}
    + \frac{1}{2}\h^{\a}\h^{\b}Q_{\b\a}^{\g}\der{}{\h^{\g}}\\+
    \left(\dot\h^{\a}Q_{\a}^a+\h^{\a}\dot Q_{\a}^a\right)\der{}{\dot
    x^a}+
    \left(\dot\h^{\a}\h^{\b}Q_{\b\a}^{\g}+
    \frac{1}{2}\h^{\a}\h^{\b}\dot Q_{\b\a}^{\g}\right)\der{}{\dot\h^{\g}}
    \,,
\end{multline}
which is the homological vector field defining the tangent
prolongation Lie algebroid in~\eqref{eq.datotbtm}. Here, expressions
such as $\dot Q_{\a}^a$ stand for $\dot x^b \p_b Q_{\a}^a$, etc.

In order to simplify notation, below we shall write  formulas for
the case when  $D, A, B, M$ are ordinary manifolds, not
supermanifolds. This will allow us to avoid some signs. Obviously,
everything carries over to the general super case.

Recall the formula for the anchor map $a\co D\to TB$.
From~\eqref{eq.qdb} we can extract
\begin{equation}\label{eq.api}
    \dot x^a=u^iQ_i^a, \quad
    \dot \h^{\b}=u^i \h^{\a}Q_{\a i}^{\b}+ \t^{\mu}Q_{\mu}^{\b}
\end{equation}
for the corresponding map $\Pi_A D\to T(\Pi B)$.

We see that the condition that~\eqref{eq.datotbtm} is a morphism of
Lie algebroids translates into the condition that the map given
by~\eqref{eq.api} relates the vector fields given by~\eqref{eq.qda}
and~\eqref{eq.qbmhat}. After a simplification, this gives the
following four equations:
\begin{equation}\label{eq.anchor1}
    Q_{\mu}^{\b}\,Q_{\b}^a=Q_{\mu}^j\,Q_j^a\,,
\end{equation}
\begin{equation}\label{eq.anchor2}
    Q_{\a i}^{\b}\,Q_{\b}^a+Q_i^b\,\p_b
    Q_{\a}^a=Q_{\a}^b\,\p_bQ_i^a+Q_{i\a}^j\,Q_j^a\,,
\end{equation}
\begin{multline}\label{eq.anchor3}
    Q_{[\a
    i}^{\d}\,Q_{\b]\d}^{\g}+\frac{1}{2}\,Q_i^b\,\p_bQ_{\b\a}^{\g}=
    Q_{[\a}^a\,\p_aQ_{\b]i}^{\g}+Q_{i[\a}^j\,Q_{\b]j}^{\g} +
\frac{1}{2}\,Q_{\b\a}^{\d}\,Q_{\d i}^{\g} \\+
\frac{1}{2}\,Q_{i\b\a}^{{\lam}}\,Q_{{\lam}}^{\g}\,,
\end{multline}
and
\begin{equation}\label{eq.anchor4}
    Q_{\mu}^{\a}\,Q_{\b\a}^{\g}=-Q_{\b}^a\,\p_aQ_{\mu}^{\g}+Q_{\mu}^j\,Q_{\b
    j}^{\g}-Q_{\mu\b}^{{\lam}}\,Q_{{\lam}}^{\g}
\end{equation}
(square brackets denote alternation).

To complete the analysis of Condition II, we have to consider
diagrams similar to~\eqref{eq.datotbtm1} and \eqref{eq.datotbtm}:
\begin{equation*}
    \begin{CD} D@>>> B\\
                @V{a}VV  @VV{a}V \\
                 TA@>>>TM
    \end{CD}
    \text{\quad and \quad}
    \begin{CD}
    \begin{CD} D@>>> B
    \end{CD}\\
    @VVV \\
    \begin{CD}  TA @>>>  TM
    \end{CD}
    \end{CD}
\end{equation*}
where $A$ and $B$ have exchanged places. This adds two equations to
the system~\eqref{eq.anchor1}--\eqref{eq.anchor4}:
\begin{multline}\label{eq.anchor5}
    Q_{[i\a}^l\,Q_{j]l}^k+\frac{1}{2}\,Q_{\a}^b\,\p_bQ_{ji}^k=Q_{[i}^a\,\p_aQ_{j]\a}^k+
    Q_{\a[i}^{\b}\,Q_{j]\b}^k+\frac{1}{2}\,Q_{ji}^l\,Q_{l\a}^k \\+\frac{1}{2}\,Q_{\a
    ji}^{{\lam}}\,Q_{{\lam}}^k
\end{multline}
and
\begin{equation}\label{eq.anchor6}
    Q_{\mu}^i\,Q_{ji}^k=-Q_j^a\,\p_aQ_{\mu}^k+Q_{\mu}^{\b}\,Q_{j\b}^k-Q_{\mu
    j}^{{\lam}}\,Q_{{\lam}}^k\,.
\end{equation}
(Note that equations~\eqref{eq.anchor1}, \eqref{eq.anchor2} are
symmetric w.r.t. $A$ and $B$.)

\textit{The system of
equations~\eqref{eq.anchor1}--\eqref{eq.anchor6} is equivalent to
Condition II of Definition~\ref{def.dlie}.} Notice that it is
bilinear in vector fields $Q_{DA}$ and $Q_{DB}$.

\subsection{Condition III, and  conclusion of the
proof}

Condition III requires ``deciphering'' more than the previous ones.
First of all, we have to consider two diagrams:
\begin{equation}\label{eq.dualda}
    \begin{CD} D^{*A}@>>> K^*\\
                @VVV  @VVV \\
                A@>>>M
    \end{CD}
\end{equation}
and
\begin{equation}\label{eq.dualdb}
    \begin{CD} D^{*B}@>>> B\\
                @VVV  @VVV \\
                K^*@>>>M
    \end{CD}
\end{equation}
and understand why the top horizontal arrow in~\eqref{eq.dualda} and
the left vertical arrow in~\eqref{eq.dualdb}  are   Lie algebroids.

Recall that the core $K$ is a vector bundle over the  base $M$. In
coordinates, $K$ is described by $x^a,z^{\mu}$ with the
transformation law
\begin{equation*}
   z^{\mu}=z^{\mu'}T_{\mu'}{}^{\mu}(x')
\end{equation*}
obtained from~\eqref{eq.lawforz} by setting $u^i$ and $w^{\a}$ to
zero. When we dualize the vertical bundle $D\to A$ in~\eqref{eq.dvb}
we obtain the bundle $D^{*A}\to A$ with fiber coordinates $w_{\a},
z_{\mu}$ (with lower indices) so that the form
$w^{a}w_{\a}+z^{\mu}z_{\mu}$ giving the pairing in coordinates, is
invariant. We arrive at the following transformation laws
\begin{align}
    w_{\a'}&=T_{\a'}^{\a} w_{\a}+u^{i'}T_{\a' i'}^{\mu}z_{\mu}\\
    z_{\mu'}&=T_{\mu'}^{\mu}z_{\mu}
\end{align}
where $(T_{\mu'}^{\mu})$ and $(T_{\mu}^{\mu'})$ are  reciprocal
matrices, and the transformation of $u^i$  remains as
in~\eqref{eq.lawforu}. This explains the double vector bundle
structure of~\eqref{eq.dualda}, in particular the vector bundle
$D^{*A}\to K^*$ (note that $z_{\mu}$ can be considered as fiber
coordinates for $K^*\to M$).

The same holds when we dualize over $B$. The total space of the
vector bundle $D^{*B}\to B$ has coordinates $x^a, u_i, w^{\a},
z_{\mu}$, the coordinates $(u_i, z_{\mu})$ being dual to $(u^i,
z^{\mu})$ on $D$. Hence the transformation law
\begin{align}
    u_{i'}&=T_{i'}^{i} u_{i}+w^{\a'}T_{\a' i'}^{\mu}z_{\mu}\\
    z_{\mu'}&=T_{\mu'}^{\mu}z_{\mu}\,,
\end{align}
from which we immediately obtain the double vector bundle structure
of~\eqref{eq.dualdb}. Treating $D^{*A}$ and $D^{*B}$ as bundles over
$K^*$, with fiber coordinates $(u^{i}, w_{\a})$ and $(u_{i},
w^{\a})$, respectively, we arrive at a surprising natural duality
between them~\cite{mackenzie:sympldouble} (see
also~\cite{mackenzie:matched} and
\cite[$\S$9.2]{mackenzie:book2005}), with the pairing   given by the
form
\begin{equation}\label{eq.pairing}
    u^iu_i-w^{\a}w_{\a}
\end{equation}
(where the minus sign is absolutely essential, in order to cancel
terms with $z_{\mu}$ appearing in  changes of coordinates).

The Lie algebroid structure  in both $D^{*A}\to K^*$ and $D^{*B}\to
K^*$ is a consequence of two facts: the above duality between
$D^{*A}\to K^*$ and $D^{*B}\to K^*$ and the linearity over the base
$K^*$ of the Poisson brackets induced on each  $D^{*A}$ and $D^{*B}$
by the Lie algebroid structures in $D\to A$ and $D\to B$,
respectively. After dualizing over $K^*$, these linear Poisson
structures give Lie algebroid structures.

All this is readily expressed in our language.

For example, the Lie algebroid structure in $D^{*B}\to K^*$ is given
by a homological vector field on $\Pi_{K^*} D^{*B}$, which is
essentially the `transpose' of the vector field $Q_{DA}$ on $\Pi_A
D$. Indeed, $Q_{DA}$  generates linear transformations of the fibers
of $\Pi_B D \to \Pi A$ and therefore   induces a vector field on the
vertical dual with reversed parity, $\Pi_{K^*} D^{*A}$. In
coordinates $x^a, \x_i, \h^{\a}, z_{\mu}$ on $\Pi_{K^*} D^{*A}$ we
get the field
\begin{multline}\label{eq.qonpikstardstara}
    Q=\h^{\a}Q_{\a}^a\der{}{x^a}-
    \left(\h^{\a}Q_{i\a}^{j}\x_j+
    \frac{1}{2}\h^{\a}\h^{\b}u^{i}Q_{i \b\a}^{\lambda}z_{{\lam}}\right)\der{}{\x_{i}}
    + \frac{1}{2}\h^{\a}\h^{\b}Q_{\b\a}^{\g}\der{}{\h^{\g}}\\ -
    \left(Q_{\mu}^{j}\x_j+\h^{\a}\t^{\mu}Q_{\mu\a}^{\lambda}z_{{\lam}}\right)
    \der{}{z_{\mu}}\,.
\end{multline}
The Schouten bracket on $\Pi_{K^*} D^{*B}$ corresponding to the Lie
algebroid structure of $D^{*A}\to K^*$ (but initially arising from
the Lie algebroid structure in $D\to B$) is defined by
\begin{equation}\label{eq.schoutenpikstardstarb}
\begin{aligned}
    \{x^a,\x_i\}&=Q_i^a\\
    \{\x_i,\x_j\}&=Q_{ij}^l z_l+\h^{\a}Q_{\a ij}^{{\lam}}z_{{\lam}}\\
    \{\x_i,\h^{\a}\}&=\h^{\b}Q_{\b i}^{\a}\\
    \{\x_i,z_{\mu}\}&=-Q_{\mu i}^{{\lam}}z_{{\lam}}\\
    \{\h^{\a},z_{\mu}\}&=-Q_{\mu}^{\a}\,,\\
\end{aligned}
\end{equation}
the other brackets between coordinates being zero.
(Formulas~\eqref{eq.qonpikstardstara} and
\eqref{eq.schoutenpikstardstarb} are given for the setup where the
initial vector bundles are purely even; the general case is similar
but contains extra signs.)

Similarly, we can write down the homological vector field and
Schouten bracket on $\Pi_{K^*} D^{*A}$. The expressions will be
`dual' to \eqref{eq.qonpikstardstara} and
\eqref{eq.schoutenpikstardstarb}.

One can immediately see that the field given
by~\eqref{eq.qonpikstardstara} is related with $Q_{BM}$ by the
projection $\Pi_{K^*} D^{*B}\to \Pi B$, which is the functor
${{\Pi}}$ applied to  the projection $D^{*B}\to B$. Hence the
horizontal projections in~\eqref{eq.dualdb} give a Lie algebroid
morphism. The same is true for the vertical projections
in~\eqref{eq.dualda}. \textit{We see that this part of Condition III
holds automatically!}

Let us examine the other part of the condition, that is,  that the
dual bundles $D^{*A}\to K^*$ and $D^{*B}\to K^*$ form a Lie
bialgebroid. It turns out to be the main one. By the self-duality of
the notion of a Lie bialgebroid, it suffices to check the two
structures on one of the bundles. Taking $D^{*B}\to K^*$, we see
that  the bialgebroid condition is equivalent to the vector field
\eqref{eq.qonpikstardstara} being a derivation of the Schouten
bracket given by~\eqref{eq.schoutenpikstardstarb}. (We have to check
not only~\eqref{eq.schoutenpikstardstarb}, but also all the zero
brackets such as $\{x^a,x^b\}=0$, etc.)

A direct calculation leads to the following set of nine equations:
\begin{equation}\label{eq.bialg1}
    Q_{\a}^a\,Q_{\mu}^{\a}-Q_{\mu}^{i}\,Q_{i}^a=0
\end{equation}
\begin{equation}\label{eq.bialg2}
    Q_{\mu}^{i}\,Q_{\nu i}^{{\lam}}+Q_{\mu}^{\a}\,Q_{\nu\a}^{{\lam}}+
    Q_{\nu}^i\,Q_{\mu i}^{{\lam}}+Q_{\nu}^{\a}\,Q_{\mu\a}^{{\lam}}=0
\end{equation}
\begin{equation}\label{eq.bialg3}
    Q_{\a
    j}^{\b}\,Q_{\b}^a+Q_j^b\,\p_bQ_{\a}^a-Q_{j\a}^i\,Q_i^a=-Q_{\a}^b\,\p_bQ_i^a
\end{equation}
\begin{equation}\label{eq.bialg4}
    Q_{j}^{a}\,\p_aQ_{\mu}^i+Q_{\mu}^k\,Q_{jk}^i-Q_{\mu}^{\a}\,Q_{j\a}^i
    =Q_{\mu j}^{{\lam}}\,Q_{{\lam}}^i
\end{equation}
\begin{multline}\label{eq.bialg5}
    -Q_{\mu}^i\,Q_{\b ji}^{{\lam}}-Q_{\b j}^{\a}\,Q_{\mu\a}^{{\lam}}+Q_{\nu
    j}^{{\lam}}\,Q_{\mu \b}^{\nu}+Q_j^a\,\p_aQ_{\mu\b}^{{\lam}}+
    Q_{\mu i}^{{\lam}}\,Q_{j\b}^i-Q_{\mu}^{\a}\,Q_{j\b\a}^{{\lam}}
    \\=
    Q_{\b}^a\,\p_aQ_{\mu j}^{{\lam}}+Q_{\mu j}^{\nu}\,Q_{\nu\b}^{{\lam}}
\end{multline}
\begin{equation}\label{eq.bialg6}
    Q_{\mu}^i\,Q_{\a
    i}^{\g}-Q_{{\lam}}^{\g}\,Q_{\mu\a}^{{\lam}}-Q_{\mu}^{\b}\,Q_{\a\b}^{\g}=
    Q_{\a}^a\,\p_aQ_{\mu}^{\g}
\end{equation}
\begin{multline}\label{eq.bialg7}
    -Q_{\a j}^{\b}\,Q_{i\b}^k-Q_{j
    l}^k\,Q_{i\a}^l-Q_{j}^a\,\p_aQ_{i\a}^k+
    Q_{\a i}^{\b}\,Q_{j\b}^k+Q_{i
    l}^k\,Q_{j\a}^l+Q_{i}^a\,\p_aQ_{j\a}^k\\
    =
    Q_{\a}^a\,\p_aQ_{ij}^k+Q_{ij}^l\,Q_{l\a}^k-Q_{\a
    ij}^{\mu}\,Q_{\mu}^k
\end{multline}
\begin{multline}\label{eq.bialg8}
    Q_{ij}^l\,Q_{l\b\a}^{{\lam}}+Q_{\b\a}^{\g}\,Q_{\g
    ij}^{{\lam}}+Q_{\a}^a\,\p_aQ_{\b ij}^{{\lam}}-Q_{\b}^a\,\p_a Q_{\a
    ij}^{{\lam}}-Q_{\a ij}^{\mu}+Q_{\b ij}^{\mu}\,Q_{\mu\a}^{{\lam}}\\
    =
    Q_{\b jl}^{{\lam}}\,Q_{i\a}^l-Q_{\a jl}^{{\lam}}\,Q_{i\b}^{{\lam}}-
    Q_{\a j}^{\g}\,Q_{i\b\g}^{{\lam}}+Q_{\b j}^{\g}\,Q_{i\a\g}^{{\lam}} -
    Q_j^a\,\p_aQ_{i\b\a}^{{\lam}}\\
    -Q_{i\b\a}^{\mu}\,Q_{\mu j}^{{\lam}}
     -Q_{\b il}^{{\lam}}\,Q_{j\a}^l+Q_{\a il}^{{\lam}}\,Q_{j\b}^{{\lam}}+
    Q_{\a i}^{\g}\,Q_{j\b\g}^{{\lam}} \\-Q_{\b i}^{\g}\,Q_{j\a\g}^{{\lam}}
    +
    Q_i^a\,\p_aQ_{j\b\a}^{{\lam}}+Q_{j\b\a}^{\mu}\,Q_{\mu i}^{{\lam}}
\end{multline}
\begin{multline}\label{eq.bialg9}
    Q_{\b k}^{\g}\,Q_{j\a}^k-Q_{\a
    k}^{\g}\,Q_{j\b}^k+Q_{j\b\a}^{{\lam}}\,Q_{{\lam}}^{\g}+
    Q_{\a j}^{\e}\,Q_{\b\e}^{\g}-Q_{\b
    j}^{\e}\,Q_{\a\e}^{\g}-Q_j^a\,\p_aQ_{\b\a}^{\g}\\
    =
    Q_{\b\a}^{\g}\,Q_{\g j}^{\b}+
    Q_{\a}^a\,\p_aQ_{\b j}^{\g}-Q_{\b}^a\,\p_aQ_{\a j}^{\g}
\end{multline}

\textit{Condition III of Definition~\ref{def.dlie} is equivalent to
equations~\eqref{eq.bialg1}--\eqref{eq.bialg9}.} Notice that they
contain equations~\eqref{eq.anchor1}--\eqref{eq.anchor6}, which are
equivalent to Condition II. Hence,  \textit{Condition III contains
Condition II.}

\smallskip
To conclude the proof of Theorem~\ref{thm.main} it remains to
consider the vector fields $Q_1$ and $Q_2$ on $\Pi^2 D$, given
by~\eqref{eq.fieldqapi} and \eqref{eq.fieldqbpi},  and  find their
commutator. A direct calculation shows that the commutativity
condition
$$[Q_1,Q_2]=0$$
produces a system of equations that coincides with
equations~\eqref{eq.bialg1}--\eqref{eq.bialg9} of Condition III (and
containing, as we showed, Condition II).

Hence we see that Mackenzie's Definition~\ref{def.dlie} is
equivalent to the commutativity of the homological fields $Q_1$ and
$Q_2$, \textbf{quod erat demonstrandum}.

\section{The big picture}\label{sec.big}

After presenting a `computational' proof of our main statement, we
shall now explain why it all works. The argument in this section
gives an alternative proof of Theorem~\ref{thm.main}, almost without
calculations.

Let us again consider a double vector bundle
\begin{equation} \label{eq.dv1}
    \begin{CD} D@>>> B\\
                @VVV  @VVV \\
                A@>>>M
    \end{CD}
\end{equation}
We shall assume that all sides of it have Lie algebroid structures
and that each of these structures is compatible with the linear
structure  in the other direction. Thus we assume the obvious part
of the definition of a double Lie algebroid (i.e., Condition I). As
we have seen, this is equivalent to saying that the Lie algebroid
structures on $D\to A$ and $D\to B$ are defined by homological
vector fields of weights $(0,1)$ and $(1,0)$ on the ultimate total
spaces of
\begin{equation} \label{eq.dv2and3}
    \begin{CD} \Pi_AD@>>> \Pi B\\
                @VVV  @VVV \\
                A@>>>M
    \end{CD}
\text{\quad and \quad}
    \begin{CD} \Pi_BD@>>>  B\\
                @VVV  @VVV \\
                \Pi A@>>>M
    \end{CD}
\end{equation}
respectively.

We want to formulate a `compatibility condition'.

Let us return for a moment to ordinary Lie algebroids (or just Lie
algebras). Suppose $E\to M$ is a vector bundle. It has three
\textit{neighbors}: the dual bundle $E^*$, the opposite bundle $\Pi
E$ and the antidual $\Pi E^*$. A Lie algebroid structure in $E$
(which is a structure on the module of sections) is equivalently
expressed by each of the following structures on its neighbors: a
homological vector field of weight $1$ on $\Pi E$, a linear Poisson
bracket on $E^*$, and a linear Schouten bracket on $\Pi E^*$. All
axioms of a Lie algebroid are contained in the equation $Q^2=0$ or
in the Jacobi identities for the Poisson or Schouten bracket.

Acting in a similar way, let us consider all the neighbors of our
double vector bundle~\eqref{eq.dv1}. There are four operations that
can be applied:  vertical dual,  horizontal dual, vertical reversion
of parity, and horizontal reversion of parity.
Besides~\eqref{eq.dv2and3} one can obtain the following double
vector bundles that are neighbors of~\eqref{eq.dv1}.

The complete parity reversion of~\eqref{eq.dv1}:
\begin{equation} \label{eq.pikvadratt}
    \begin{CD} \Pi^2D@>>>  \Pi B\\
                @VVV  @VVV \\
                \Pi A@>>>M
    \end{CD}
\end{equation}
The two duals of~\eqref{eq.dv1}:
\begin{equation} \label{eq.dv5and6}
    \begin{CD} D^{*A}@>>>  K^*\\
                @VVV  @VVV \\
                A@>>>M
    \end{CD}
\text{\quad and \quad}
    \begin{CD} D^{*B}@>>>  B\\
                @VVV  @VVV \\
                K^*@>>>M
    \end{CD}
\end{equation}
The parity reversions of each of the duals:
\begin{equation} \label{eq.dv7}
    \begin{CD} \Pi_AD^{*A}@>>>  \Pi K^*\\
                @VVV  @VVV \\
                A@>>>M
    \end{CD}
\end{equation}
\begin{equation} \label{eq.piduala}
    \begin{CD} \Pi_{K^*}D^{*A}@>>>   K^*\\
                @VVV  @VVV \\
                \Pi A@>>>M
    \end{CD}
\end{equation}
\\
\begin{equation} \label{eq.pikvadratduala}
    \begin{CD} \Pi^2 D^{*A}@>>>   \Pi K^*\\
                @VVV  @VVV \\
                \Pi A@>>>M
    \end{CD}
\end{equation}
and
\begin{equation} \label{eq.pidualb}
    \begin{CD} \Pi_{K^*}D^{*B}@>>>  \Pi B\\
                @VVV  @VVV \\
                K^*@>>>M
    \end{CD}
\end{equation}
\\
\begin{equation} \label{eq.dv11}
    \begin{CD} \Pi_{B}D^{*B}@>>>   B\\
                @VVV  @VVV \\
               \Pi K^*@>>>M
    \end{CD}
\end{equation}
\\
\begin{equation} \label{eq.pikvadratdualb}
    \begin{CD} \Pi^2 D^{*B}@>>>   \Pi B\\
                @VVV  @VVV \\
                \Pi K^*@>>>M
    \end{CD}
\end{equation}
This is the full list up to natural isomorphisms. These twelve
objects, including the original double vector bundle~\eqref{eq.dv1},
can be arranged into a four-valent graph. More precisely, from each
vertex emanate two edges corresponding to taking duals and two edges
corresponding to parity reversions. We can think of the edges being
colored by two colors.

\begin{rem} In the multiple case, for an $n$-fold vector bundle, the
number of edges emanating from each vertex will be $n+n=2n$.
Question: what is the total number of neighbors (the number of
vertices in the graph)?
\end{rem}

Lie algebroid structures on the sides of ~\eqref{eq.dv1}  satisfying
the linearity conditions, which  were expressed above in terms of
weights, generate various combinations of structures on each of the
neighbors~\eqref{eq.dv2and3}--\eqref{eq.pikvadratdualb}, a pair for
each. More symmetrically, we can say that each pair of structures
for a particular double vector bundle
from~\eqref{eq.dv1}--\eqref{eq.pikvadratdualb} is just a
manifestation of one `double structure'. (All pairs contain the same
information.) One can make a list of such structures. The next step
will be to look for suitable compatibility conditions for each pair.
The philosophy is  that one should look for pairs where a
compatibility condition is formulated naturally, and take it as the
definition of  compatibility for an equivalent pair where such a
condition does not come about in an obvious way. Suppose we do not
know what  a double Lie algebroid is (a compatibility of Lie
algebroids on the sides of~\eqref{eq.dv1}): for the right notion,
examine its neighbors.

For each of the double bundles~\eqref{eq.dv2and3}, the ultimate
total space is a Lie algebroid (over the base $K^*$) and
simultaneously possesses a linear Poisson bracket (linear over both
bases). That means that the dual bundle over $K^*$ is also a Lie
algebroid, and one may ask whether they form a Lie bialgebroid. As
the analysis of the previous section shows, this may be considered
as the \textbf{Mackenzie definition} of a double Lie algebroid
for~\eqref{eq.dv1}, since it subsumes all other conditions from
Definition~\ref{def.dlie}.

As the example of Lie algebras (or algebroids) shows, to define a
Lie algebra or algebroid in terms of a structure on the manifold
itself~\footnote{We mean a structure in the algebra of functions as
opposed to, say, a structure on sections of a vector bundle.}, one
has to replace $E$ by one of the neighbors: $E^*$, $\Pi E$, or $\Pi
E^*$, with the corresponding structure. For the ``bi-'' case (Lie
bialgebras or bialgebroids),  remarkably, it all boils down to one
type of structure, should we take $E^*$, $\Pi E$, or $\Pi E^*$ as a
model: namely, a $QS$-structure, i.e., a homological vector field
and a Schouten bracket, with the derivation condition
(see~\cite{tv:graded}).

Among the twelve double vector
bundles~\eqref{eq.dv1}--\eqref{eq.pikvadratdualb} there are
precisely five with a structure induced on the  total space, thus
allowing a verifiable compatibility condition. They
are:~\eqref{eq.pikvadratt}, \eqref{eq.piduala},
\eqref{eq.pikvadratduala}, \eqref{eq.pidualb},
\eqref{eq.pikvadratdualb}.

On the total space of~\eqref{eq.pikvadratt} there are two
homological vector fields of weights $(0,1)$ and $(1,0)$. A
compatibility condition for them is commutativity. This is precisely
our Definition~\ref{def.dalie}.

On the total space  of~\eqref{eq.piduala} or \eqref{eq.pidualb}
there is a Schouten bracket of weight $(-1,-1)$ and a homological
vector field of weight $(0,1)$ or $(1,0)$, respectively. On the
total space of~\eqref{eq.pikvadratduala} or
\eqref{eq.pikvadratdualb} there is a Poisson bracket  of weight
$(-1,-1)$ and a homological vector field of weight $(0,1)$ or
$(1,0)$. A compatibility condition in each case is the derivation
property of the vector field w.r.t. the bracket.

\begin{prop} \label{prop.bialg}
The compatibility conditions for \eqref{eq.piduala},
\eqref{eq.pidualb}, \eqref{eq.pikvadratduala} and
\eqref{eq.pikvadratdualb} are equivalent, and are different ways of
saying that $(D^{*A}, D^{*B})$ is a Lie bialgebroid over $K^*$.
\end{prop}
\begin{proof} Indeed, for any Lie bialgebroid $(E,E^*)$ the compatibility can
be stated in terms of either $E$ or $E^*$ (as a $QS$-structure on
either $\Pi E$ or $\Pi E^*$, respectively). This corresponds
to~\eqref{eq.piduala} or \eqref{eq.pidualb} in our case. In our
special case  there is also an extra option of changing parity in
the second direction, which adds~\eqref{eq.pikvadratduala} and
\eqref{eq.pikvadratdualb} to the picture.
\end{proof}

We see now that there are essentially two conditions to compare: the
Mackenzie bialgebroid condition, which lives on  one
of~\eqref{eq.piduala}, \eqref{eq.pidualb}, \eqref{eq.pikvadratduala}
and \eqref{eq.pikvadratdualb}, and our commutativity condition
for~\eqref{eq.pikvadratt}.

\begin{prop} \label{prop.equiv}
The Mackenzie bialgebroid condition and the commutativity of the two
homological vector fields  for~\eqref{eq.pikvadratt} are equivalent.
\end{prop}
\begin{proof} Consider one of the manifestations of the bialgebroid
condition, say, for concreteness,~\eqref{eq.piduala}. The derivation
property means that the flow of the vector field preserves the
bracket. On the other hand, the commutativity condition
for~\eqref{eq.pikvadratt} means that the flow of one field preserves
the other. Now the claim follows from  functoriality: notice that a
linear transformation preserves a Lie bracket if and only if the
adjoint map preserves the corresponding linear Poisson bracket and
if and only if the `$\Pi$-symmetric' map preserves the corresponding
homological vector field.
\end{proof}

Propositions~\ref{prop.bialg} and~\ref{prop.equiv} together imply
Theorem~\ref{thm.main}.

\section{Applications and generalizations}\label{sec.appl}

In this section we consider some examples and applications, as well
as discuss an extension to multiple Lie algebroids.

\subsection{Doubles of Lie bialgebroids.}

Recall that Drinfeld's \textit{classical double} of a Lie bialgebra
is again a Lie bialgebra with ``good'' properties. An analog of this
construction for Lie algebroids turned out to be a puzzle.  Three
constructions of a `double' of a Lie bialgebroid have been
suggested. Suppose $(E,E^*)$ is a Lie bialgebroid over a base $M$.
Liu, Weinstein and Xu~\cite{weinstein:liuxu} suggested to consider
as its double a structure of a Courant algebroid on the direct sum
$E\oplus E^*$. Mackenzie in~\cite{mackenzie:doublealg,
mackenzie:doublealg2, mackenzie:drinfeld, mackenzie:notions} and
Roytenberg in~\cite{roytenberg:thesis} suggested two different
constructions based on cotangent bundles. Though they look very
different (in particular, Roytenberg's double is a supermanifold,
and Mackenzie stays in the classical world), we shall  show now that
they are essentially the same.

Roytenberg previously showed~\cite{roytenberg:thesis} that the
Liu--Weinstein--Xu double is recovered from his own construction as
a derived bracket, generalizing the results  of
C.~Roger~\cite{roger:1991} and
Y.~Kosmann-Schwarzbach~\cite{yvette:jacobian, yvette:derived} for
Lie bialgebras. Therefore, proving that the Mackenzie and Roytenberg
pictures  are equivalent or, actually,   the same, if understood
properly, shows conclusively that this `cotangent double' is
fundamental, and should be regarded as the correct extension of
Drinfeld's double of Lie bialgebras to Lie bialgebroids.

Both Roytenberg's and Mackenzie's construction use the statement
that the cotangent bundles of dual vector bundles are isomorphic
(\cite{mackenzie:bialg}, an extension of ~\cite{tulczyjew:1977}; see
also~\cite{mackenzie:diffeomorphisms}, \cite{roytenberg:thesis},
\cite{tv:graded}). Hence there is a double vector bundle
\begin{equation} \label{eq.mackdouble}
    \begin{CD} T^*E=T^*E^*@>>>   E^*\\
                @VVV  @VVV \\
               E@>>>M
    \end{CD}
\end{equation}
Mackenzie shows that it is a double Lie algebroid. He calls it the
\textit{cotangent double} of a Lie bialgebroid $(E,E^*)$. Note that
the canonical symplectic structure on $T^*E$ corresponds to the
invariant scalar product on Drinfeld's double $\mathfrak d
(\mathfrak b)=\mathfrak b\oplus \mathfrak b^*$ of a Lie bialgebra
$\mathfrak b$.

On the other hand, Roytenberg uses the description of Lie algebroids
via homological vector fields. He considers the double vector bundle
\begin{equation} \label{eq.roytdouble}
    \begin{CD} T^*\Pi E=T^*\Pi E^*@>>>   \Pi E^*\\
                @VVV  @VVV \\
               \Pi E@>>>M
    \end{CD}
\end{equation}
and   homological vector fields  $Q_E\in\Vect (\Pi E)$ and
$Q_{E^*}\in\Vect (\Pi E^*)$ defining Lie algebroid structures on
$E\to M$ and $E^*\to M$, respectively. Recall that  vector fields on
a manifold correspond to fiberwise linear functions (Hamiltonians)
on the cotangent bundle so that the commutator  maps to the Poisson
bracket. Denote the functions corresponding to $Q_E$ and $Q_{E^*}$
by $H_E$ and $H_{E^*}$, respectively. Roytenberg shows that under
the natural symplectomorphism $T^*\Pi E=T^*\Pi E^*$ the linear
function $H_{E^*}$ on $T^*\Pi E^*$ corresponding to the vector field
$Q_{E^*}$ transforms precisely into the fiberwise quadratic function
$S_E$ on $T^*\Pi E$  specifying the Schouten bracket on $\Pi E$
induced by the Lie structure on $E^*$. The derivation property of
$Q_E$ w.r.t. the Schouten bracket on $\Pi E$ is one of the
equivalent definitions of a Lie bialgebroid~\cite{yvette:exact}, and
the most convenient. Hence, Roytenberg's statement means that it is
also equivalent to the commutativity of the Hamiltonians $H_E$ and
$H_{E^*}$ under the canonical Poisson bracket. They generate
commuting homological vector fields $X_{H_E}$ and $X_{H_{E^*}}$ on
the cotangent bundle $T^*\Pi E$. In our language, $X_{H_E}$ and
$X_{H_{E^*}}$ make~\eqref{eq.roytdouble} a Lie antialgebroid. (One
can see that the conditions for weights are satisfied.)
Roytenberg~\cite{roytenberg:thesis} calls the supermanifold $T^*\Pi
E=T^*\Pi E^*$ together with the homological vector field
$Q=X_{H_E}+X_{H_{E^*}}$ on it, the \textit{Drinfeld double} of
$(E,E^*)$.

If we slightly refine Roytenberg's picture,  considering the double
Lie antialgebroid given by $X_{H_E}$ and $X_{H_{E^*}}$ rather than a
single $Q$-manifold, we can immediately see that by
Theorem~\ref{thm.main} his picture becomes identical to that of
Mackenzie.

Indeed, apply the complete reversion of parity
to~\eqref{eq.mackdouble}. Notice that $\Pi^2 T^*E=\Pi^2 T^*E^*$
coincides with $T^* \Pi E=T^* \Pi E^*$ (easily checked in
coordinates). By Theorem~\ref{thm.main}, the double vector
bundle~\eqref{eq.mackdouble} is a double Lie algebroid if and only
if the corresponding double vector bundle
\begin{equation} \label{eq.mackdoublepikv}
    \begin{CD} \Pi^2T^*E=\Pi^2T^*E^*@>>>   \Pi E^*\\
                @VVV  @VVV \\
               \Pi E@>>>M
    \end{CD}
\end{equation}
which is identical with~\eqref{eq.roytdouble}, is a double Lie
antialgebroid. It remains to identify the respective homological
vector fields on the ultimate total space, which can be achieved by
a direct inspection.

We have arrived at the following statement.
\begin{prop} Roytenberg's and Mackenzie's pictures give the same notion of a double of a Lie
bialgebroid (up to a change of parity).
\end{prop}

We can now identify the two constructions and speak simply of the
(cotangent) \textit{double} of a Lie bialgebroid  as a double Lie
algebroid, most efficiently described in the ``anti-'' language of
diagram ~\eqref{eq.roytdouble}.

\subsection{More on doubles.}

Recall that Drinfeld's classical double of a Lie bialgebra is not
just a Lie algebra, but also a coalgebra, and in fact a Lie
bialgebra again. This gives a direction in which to look in the case
of Lie bialgebroids. Note that this second structure (for doubles of
Lie bialgebroids) has not been discovered previously.

However, for many people including the author there was absolutely
no doubt that such a structure exists. The following conjectured
statement was put down in my notes of 2002 as a guideline for (not
yet obtained at that time) alternative description of double Lie
algebroids.

\begin{genprinc}
Taking the double of an $n$-fold Lie bialgebroid should give  an
$(n+1)$-fold Lie bialgebroid, with an additional property, such as a
symplectic structure.
\end{genprinc}

Of course it involved new notions yet to be defined. With the
efficient description of double Lie algebroids in terms of
supermanifolds obtained in this paper, it becomes possible. Multiple
Lie algebroids can now be easily introduced  by extending our
description of the double case. See Remark~\ref{rem.mantialg} above
and the next subsection. As for ``double Lie \textbf{bi}algebroids''
(or ``bi- double Lie algebroids''), and the further multiple ``bi-''
case, this  notion is  properly defined in our joint work with
Kirill Mackenzie, and is a subject of our forthcoming
paper~\cite{mackenzie:bidouble}, where precise definitions and
statements can be found.

\begin{ex} \label{ex.double}
Consider again the double vector bundle given
by~\eqref{eq.mackdouble}. Notice that the core of it is the
cotangent bundle $T^*M\to M$. Take the two duals
of~\eqref{eq.mackdouble}. We obtain
\begin{equation} \label{eq.mackdouble3}
    \begin{CD} TE@>>>   T M\\
                @VVV  @VVV \\
               E@>>>M
    \end{CD}
\end{equation}
for the vertical dual and
\begin{equation} \label{eq.mackdouble4}
    \begin{CD} TE^*@>>>   E^*\\
                @VVV  @VVV \\
               T M@>>>M
    \end{CD}
\end{equation}
for the horizontal dual. Both double vector
bundles~\eqref{eq.mackdouble3} and \eqref{eq.mackdouble4} are known
to be double Lie algebroids~\cite{mackenzie:notions}. Three double
Lie algebroids~\eqref{eq.mackdouble}, ~\eqref{eq.mackdouble3},
~\eqref{eq.mackdouble4}, which are double vector bundles in duality
\begin{equation}
    \begin{picture}(300,80)(0,65)
    \put(90,90){${\begin{CD} {TE}  @>>>   T M\\
                @VVV  @VVV \\
               E@>>>M\end{CD}}$}
    \put(120,110){$\begin{CD} {} @.   T E^*\\
                @.  @VVV \\
               T^*E@>>>E^*\end{CD}$}
    {\put(132,84){\vector(-2,-1){18}}}
    {\put(192,84){\vector(-2,-1){18}}}
    {\put(192,130){\vector(-2,-1){18}}}
    \end{picture}
\end{equation}
make a \textit{double bialgebroid} (or  `bi double Lie algebroid'),
for example, in the sense that the triple vector bundle
\begin{equation}
    \begin{picture}(300,80)(0,65)
    \put(90,90){${\begin{CD} {TE}  @>>>   T M\\
                @VVV  @VVV \\
               E@>>>M\end{CD}}$}
    \put(120,110){$\begin{CD} T^*T^*E @>>>   T E^*\\
                @VVV  @VVV \\
               T^*E@>>>E^*\end{CD}$}
    {\put(132,84){\vector(-2,-1){18}}}
    {\put(192,84){\vector(-2,-1){18}}}
    \put(132,130){\vector(-2,-1){18}}
    {\put(192,130){\vector(-2,-1){18}}}
    \end{picture}
\end{equation}
where they make a corner, is a triple Lie algebroid. (It is the
second cotangent double of $(E,E^*)$.)
\end{ex}

This example should be considered as a preliminary announcement of
results to follow in~\cite{mackenzie:bidouble}.

\subsection{Higher Lie algebroids}

We have already mentioned that the methods of this paper naturally
allow us  to consider multiple Lie algebroids. (A part of the
motivation for doing this comes from the theory of doubles discussed
above.) A detailed theory will be presented elsewhere. Here we wish
to give a sketch.

First we need a language for describing multiple vector bundles. Fix
a natural number $n$. To define $n$-fold vector bundles, consider
vector spaces $V_r$, $V_{r_1r_2}$, $V_{r_1r_2r_3}$, \ldots, of
arbitrary dimensions $d_r$,$d_{r_1r_2}$, $d_{r_1r_2r_3}$, etc.,
numbered by increasing sequences $r_1<\ldots <r_k$, where $0<k\leq
n$ and all $r_i$ run from $1$ to $n$.

\begin{ex} When $n=1$, we have just one vector space $V=V_1$. When
$n=2$, we have $V_1$, $V_2$ and $V_{12}$. For $n=3$, we have $7$
spaces: $V_1$, $V_2$, $V_3$, $V_{12}$, $V_{13}$, $V_{23}$, and
$V_{123}$. In general the number of spaces is $2^n-1$.
\end{ex}

For convenience of notation let us fix linear coordinates on each of
the spaces, denoting them $v^{i_r}_{(r)}$,
$v^{i_{r_1r_2}}_{(r_1r_2)}$, etc. (Each index such as $i_r$ runs
over its own set of values, of cardinality equal to the dimension of
the respective space.)

\begin{de} An \textit{$n$-fold vector bundle} over a base $M$ is a fiber
bundle over $M$ with the standard fiber
\begin{equation*}
    \prod_{r} V_r \,\times \,\prod_{r_1<r_2} V_{r_1r_2} \,\times \,\ldots
    \,\times\,
    V_{12\ldots n}
\end{equation*}
where the transition functions have the form:
\begin{align*}
    v^{i_r}_{(r)}&=v^{i_r'}_{(r)} T_{i_r'}^{i_r},  \\
    v^{i_{r_1r_2}}_{(r_1r_2)}&=v^{i_{r_1r_2}'}_{(r_1r_2)} T_{i_{r_1r_2}'}^{i_{r_1r_2}}+
    v^{i_{r_1}'}_{(r_1)}v^{i_{r_2}'}_{(r_2)}
    T_{i_{r_2}'i_{r_1}'}^{i_{r_1r_2}},
     \intertext{\quad\quad\quad\quad\quad\quad\quad\dots\dots\dots\dots}
    v^{i_{12\ldots n}}_{(12\ldots n)}&=v^{i_{12\ldots n}'}_{(12\ldots
    n)}T_{i_{12\ldots n}'}^{i_{12\ldots n}}+\ldots + v^{i_1'}_{(1)}\ldots
    v^{i_n'}_{(n)} T_{{i_n'}\ldots{i_1'}}^{i_{12\ldots n}}\,.
\end{align*}
\end{de}

In other words, the transformation for each of $V_r$ is linear; for
$V_{r_1r_2}$ it is linear plus an extra term bilinear in $V_{r_1}$
and $V_{r_2}$, etc.

\begin{ex} For a triple vector bundle ($n=3$), we have fiber coordinates:
$v^{i_1}_{(1)}$, $v^{i_2}_{(2)}$, $v^{i_3}_{(3)}$,
$v^{i_{12}}_{(12)}$, $v^{i_{13}}_{(13)}$, $v^{i_{23}}_{(23)}$, and
$v^{i_{123}}_{(123)}$. The transformation law is as follows:
\begin{align*}
    v^{i_1}_{(1)}&=v^{i_1'}_{(1)}T_{i_1'}^{i_1}\\
    v^{i_2}_{(2)}&=v^{i_2'}_{(2)}T_{i_2'}^{i_2}\\
    v^{i_3}_{(3)}&=v^{i_3'}_{(3)}T_{i_3'}^{i_3}\\
    v^{i_{12}}_{(12)}&=v^{i_{12}'}_{(12)}T_{i_{12}'}^{i_{12}}+
    v^{i_1'}_{(1)}v^{i_2'}_{(2)}T_{{i_2'}{i_1'}}^{i_{12}}\\
    v^{i_{13}}_{(13)}&=v^{i_{13}'}_{(13)}T_{i_{13}'}^{i_{13}}+
    v^{i_1'}_{(1)}v^{i_3'}_{(3)}T_{{i_3'}{i_1'}}^{i_{13}}\\
    v^{i_{23}}_{(23)}&=v^{i_{23}'}_{(23)}T_{i_{23}'}^{i_{23}}+
    v^{i_2'}_{(2)}v^{i_3'}_{(3)}T_{{i_3'}{i_2'}}^{i_{23}}\\
    v^{i_{123}}_{(123)}&=
                    \begin{aligned}[t]
    v^{i_{123}'}_{(123)}T_{i_{123}'}^{i_{123}}+
    v^{i_1'}_{(1)}v^{i_{23}'}_{(23)}T_{{i_{23}'}{i_1'}}^{i_{123}}+
    v^{i_2'}_{(2)}v^{i_{13}'}_{(13)}T_{{i_{13}'}{i_2'}}^{i_{123}}+
    v^{i_3'}_{(3)}v^{i_{12}'}_{(12)}T_{{i_{12}'}{i_3'}}^{i_{123}} \\+
v^{i_1'}_{(1)}v^{i_2'}_{(2)}v^{i_3'}_{(3)}T_{{i_3'}{i_2'}{i_1'}}^{i_{123}}
\end{aligned}
\end{align*}
\end{ex}

\begin{rem} Triple vector bundles\,---\,with the quaternary case
briefly mentioned\,---\,were introduced and studied
in~\cite{mackenzie:duality} from a different viewpoint (not using
local trivializations and transition functions).
Paper~\cite{mackenzie:duality} also contains conjectured `likely
principles' of duality for general multiple case.
\end{rem}

A multiple vector bundle has \textit{faces}, which are also multiple
vector bundles. A face is obtained by choosing indices $r_1<\ldots <
r_k$; fiber coordinates for it will be the coordinates
$v^{i_{{r_1\ldots r_k}}}_{(r_1\ldots r_k)}$ and all other
coordinates with indices labelled by subsets of $r_1,\ldots, r_k$.
For example, for a triple vector bundle there are faces that are
(ordinary) vector bundles and double vector bundles, corresponding
to the edges and $2$-faces of a $3$-cube. In a natural way  various
partial projections and zero sections are defined.

The total space of a multiple vector bundle is a multi-graded
manifold. More precisely, there are weights $\boldsymbol{w}_r$,
$r=1,\ldots, n$, each of them being a degree  in all coordinates
containing a given label $r$. For example, $\boldsymbol{w}_2$ is the
total degree in $v_{(2)}$, $v_{(12)}$, $v_{(23)}$, \ldots,
$v_{(12\ldots n)}$. We define \textit{total weight} as
$\boldsymbol{w}=\boldsymbol{w}_1+\ldots+\boldsymbol{w}_n$.

Due to the multilinearity of transition functions, for a multiple
vector bundle  the operations of \textit{partial parity reversion}
$\Pi_r$ and   \textit{partial dual} $\mathrm{D}_r$ in the $r$-th
direction, make sense for each $r=1,\ldots, n$.

\begin{de} \label{def.multalie}An \textit{$n$-fold Lie antialgebroid} is an $n$-fold
vector bundle endowed with a homological vector field $Q$  of total
weight $1$ on the total space $E$.
\end{de}
Clearly, this is the same as having $n$ odd vector fields $Q_r$ of
weights $(0,\ldots, 1, \ldots, 0)$ such that
\begin{equation*}
    [Q_r,Q_s]=0
\end{equation*}
for all $r,s$.

\begin{de} An \textit{$n$-fold Lie algebroid} is an $n$-fold
vector bundle  such that the $n$-fold vector bundle obtained by the
complete parity reversion $\Pi^n=\Pi_1\ldots \Pi_n$ is an $n$-fold
Lie antialgebroid.
\end{de}

In other words, we take the statement of Theorem~\ref{thm.main} as a
definition for the multiple case.

Each face of a multiple Lie (anti)algebroid is also a multiple Lie
(anti)algebroid.

We expect the following to be true: the possibility to define
multiple Lie algebroids \textit{\`{a}  la} Mackenzie, via duals and
bialgebroids, and the equivalence of such a definition with
Definition~\ref{def.multalie} (i.e., the analog of
Theorem~\ref{thm.main}). This requires an analysis of the structures
induced on all neighbors of a multiple Lie algebroid.

An \textit{$n$-fold Lie bialgebroid} (or a \textit{``bi-'' $n$-fold
Lie algebroid}) can be defined, in supergeometry terms, as an
$n$-fold Lie algebroid such that all its duals are also $n$-fold Lie
algebroids with the compatibility condition that reads  as follows:
on the total space with completely reversed parity there is a
homological vector field $Q$ of total weight $1$, which defines the
algebroid structure, and an odd or even (depending on the parity of
the number $n$) bracket of an appropriate weight, which corresponds
to algebroid structures on all the duals;  the field $Q$ should be a
derivation of the bracket. Shortly, it can be described as a $QS$-
or $QP$-structure on the total space with particular conditions on
weights.

It should be possible to prove that this is equivalent to the
cotangent double being an $(n+1)$-fold Lie algebroid. (Which we used
as a definition in Example~\ref{ex.double}.) The main statement then
should be that the \textbf{general principle} stated above holds:
that the cotangent double is, moreover, a ``bi-'' $(n+1)$-fold Lie
algebroid. See~\cite{mackenzie:bidouble} for all this.


\def\cprime{$'$} \def\cprime{$'$} \def\cprime{$'$}

\end{document}